\documentclass[review]{elsarticle}
\usepackage{amsthm,amssymb,amsmath,amsfonts}
\usepackage{tikz-cd} 
\usepackage{lineno,hyperref}
\modulolinenumbers[5]
\theoremstyle{definition}

\newtheorem{theorem}{Theorem}[section]
\newtheorem{corollary}{Corollary}[theorem]
\newtheorem{lemma}[theorem]{Lemma}
\newtheorem{remark}{remark}[section]
\journal{ArXiv}

\begin{document}

\begin{frontmatter}

\title{Triple extension of Tietze theorem and Baer criterion}

\author{Assad Rashidi\fnref{myfootnote}}
\address{University of Kurdistan}
\fntext[myfootnote]{a.rashidi@sci.uok.ac.ir}

\author{Kaveh Mohammadi\fnref{myfootnote}}
\address{Amirkabir University of Technology}
\fntext[myfootnote]{kmsmath@aut.ac.ir}

\begin{abstract}
In this paper, through the combination of Tietze extension theorem and Baer criteria, we build a new mathematical structure which is similar to a triangular pyramid, and then we prove that the topological space which  we call it $T_b$  appeared as a result of the combination  and sat at the apex of the pyramid is a tychonof space. Finally, we obtain three new extension theorems.
\end{abstract}


\end{frontmatter}

\section{Introduction}
Let $C$ be a closed subset of a normed space $N$ and suppose $g:C\rightarrow \mathbb{R}$ is continuous. Then Tietze's theorem asserts that g can be extended to a continuous function G defiend on all of N.\cite{andrewm.brucknerjudithb.brucknerbrians.thomson1996}
\[G|_C=g \quad \text{or}\quad Gt=g\]
\begin{center}

	\begin{tikzcd}[column sep=small]
	C\arrow{rr}{t} \arrow[swap]{dr}{g}& &N \arrow{dl}{G}\\
	& \mathbb{R} & 
	\end{tikzcd} 
\end{center}

Baer citerion states that  a left R-module E is injective if and only if any homomorphism $f : I \rightarrow E$ defined on a left ideal I of R can be extended to all of R.\cite{davids.dummitrichardm1974} \cite{johndaun1994}

\[F|_I=f \quad \text{or}\quad FT=f\]

\begin{center}
	\begin{tikzcd}
	0 \arrow[r] & I \arrow[rr, "T"] \arrow[rd, "f"] &   & \mathbb{R} \arrow[ld, "F"'] \\
	&                                   & E &                            
	\end{tikzcd}
\end{center}

By combining Tietze extension theorem and Baer criteria, we create a new mathematical structure which is similar to a triangular pyramid. There are two immediate consequences after combining these theorem. First, we prove the $T_b$ space which sits at the apex of the pyramid is a Tychonoff space. Secondly, we obtain three new extension theorems which we will discuss them in the third chapter.

\section{Statement of the idea}
Since a discrete topology can be defined under infinite algebraic structures, hence a natural topology is defined which is compatible with continuous operation of algebraic structures such as Topological groups, Topological rings, Topological vector space, etc \cite{l.gillmanm.jerison1960} \cite{sashokalajdzievski2015}. Essentially every discrete topologiccal space satisfies each of the seperation axixoms so we can assume that the topology on $T_b$ and N discrete. 

\begin{lemma}
	we can define a discrete topology on closed set.
	\begin{proof}
		As we know the union of infinitely many closed sets is not closed. Let suppose 
		$C=C_1=...=C_i$ , $i\in I$ ((I is a infinite index set))$\Longrightarrow C=\bigcup_{i=1}^{\infty} C_{i}= \bigcap^{i=1}_{\infty}C_i$, every $C_i$ is a equal copy of the closed set C and each this copies cover each other so that we can define discrete topology on C (C is infinite set of the equal copy of C itself)
		
	\end{proof}
\end{lemma}

To convey our idea, we need to define the  new notations which are necessary to use Tietze extension theorem in the combined diagram of Tietze and Baer in the following ways

\begin{itemize}
	\item $I_g=g$,  is a left ideal from R in Baer criterion which is modified to $I_g$ that is $I_g : C\rightarrow \mathbb{R}$ or $C(c)$(Left Ideal of real valued continous function from closed subset C to $\mathbb{R}$ )
	\item $R_G=G$,  is a ring in Baer criterion which is modified to $R_G$ that is $R_G: N\rightarrow \mathbb{R}$ or $C(N)$ (Ring of real valued continous function from normal space to $\mathbb{R}$)
	
\end{itemize}
Between two topological spaces which the first topological space is endowed with discrete topology, an algebraic structure can be placed as a set of continuous functions. Since E is an injective R-module we can put it between two topological spaces which one of them is $T_b$ which is endowed with discrete topology  and the other is $\mathbb{R}$ as a continuous function.\cite{willardstephen2004}\cite{l.gillmanm.jerison1960}

As we learned that left Ideal I and the Ring $\mathbb{R}$ from Baer criterion in Tietze extension Theorem can be used. The next step is to define the map 
\begin{itemize}
	\item $E:T_b\rightarrow \mathbb{R}$ or $C(T_b)$(left injective of real valued continuous function from $T_b$ to $\mathbb{R}$).
\end{itemize}
\begin{center}
	\begin{tikzcd}
	C(c) \arrow[rr, "T"] \arrow[rd, "f"'] &        & C(N) \arrow[ld, "F"] \\
	& C(T_b) &                     
	\end{tikzcd}
\end{center}

Because we can put a discrete topology on C, N, and $T_b$, then there exist continuous functions like $\alpha,\beta$, and t between these space.

\begin{center}
	\begin{tikzcd}
	C \arrow[rr, "t"] \arrow[rd, "\alpha"'] &     & N \arrow[ld, "\beta"] \\
	& T_b &                      
	\end{tikzcd}
\end{center}

After combining the Tietze extension theorem and Baer criterion, a triangular pyramid will appear in which $T_b$ is a tychonoff space. Another interesting result for this combination is that every lateral face of the triangular pyramid will appear as an extension theorem in which one of them is Tietze theorem and remaining lateral faces can be considered as three new extension theorem in the branch of topology. We will discuss these theorem in the next section.

\begin{center}

	\begin{tikzcd}
	C \arrow[rrrr, "t"] \arrow[rrd, "\alpha"] \arrow[rrdddd, "I_g" description] &                                     &                 &                     & N \arrow[lld, "\beta"'] \arrow[lldddd, "R_G" description] \\
	&                                     & T_b \arrow[ddd] &                     &                                                           \\
	O \arrow[r]                                                                 & Ig \arrow[rr, "T"'] \arrow[rd, "f"] &                 & RG \arrow[ld, "F"'] &                                                           \\
	&                                     & E               &                     &                                                           \\
	&                                     & \mathbb{R}      &                     &                                                          
	\end{tikzcd}
\end{center}
\section{Consequence of the combintation}
\begin{corollary}
	$T_b $ is a Tychonoff space.
	
	\begin{proof}
		Firstly, $T_b$ is $T_0$ because we can see distinct points in $T_b$ space which are topologically distinguishable by assuming discrete topology over it.
		
		Secondly, by the categorial  point of view, we can see the following equalizer structure of the triangular pyramid.
		
		\begin{center}

			\begin{tikzcd}
			N \arrow[r] \arrow[r, "\beta"]                 & T_b \arrow[r, "E"] & \mathbb{R} \\
			C \arrow[ru, "\alpha"'] \arrow[u, "t", dotted] &                    &  
			\end{tikzcd}
		\end{center}

		By the definition of equaliser, equaliser includs  $\forall i,j\in E$    an object N and a morphism $\beta  :N\rightarrow T_b$ ($i,j,\beta$ are continous) satisfying $io\beta=jo\beta$, and such that given any object C and morphism $\alpha:C \rightarrow T_b$, if $io\alpha=jo\alpha$ then by titze extension theorem there exist a unique morphism $t:C\rightarrow N$ such that $\beta t= \alpha$.

		According to the definition of equalizer in set theory we have $N=Eq(E)=\{\ X\in T_b |\forall i,j\in E, i(X)=j(X)\}$ and set points of X are closed subset of N (we know t is an inclusion map, and $\beta$ is a monomorphism) which are induced by $\beta $ to $T_b$, thus based on the properties of equalizer, all closed sets of $T_b$  are mapped to the point of $\mathbb{R}$ space by continuous functions of $E$. Namely, a disjoint closed subset of $T_b$ cannot be separated by a continuous function. In result, $T_b$ space is not a normal space although if given any point $y\in T_b$ such that $y\notin X$ simply there exists a continuous function of E from y to another point of $\mathbb{R}$ namely  they are separated by continuous f then $T_b$ space is completely regular.
		So $T_b$ is $T_0$ and compeletly regular space and in consequcne $T_b$ is Tychonoff space.

	\end{proof}
\end{corollary}

\subsection{Extension Theorems}

\begin{theorem}
	If N is a normal topological space and $\alpha :C \rightarrow  T_b$ is a continuous function  from C which closed subset of N to topological space $T_b$ then there exists a continuous function 
	$\beta: N\rightarrow T_b$ such that  $\beta|_C=\alpha, or \beta t
	=\alpha$
	\begin{center}

		\begin{tikzcd}
		C \arrow[rr, "t"] \arrow[rd, "\alpha"] &     & N \arrow[ld, "\beta"'] \\
		& T_b &                       
		\end{tikzcd}
	\end{center}
	
	\begin{remark}
		Based on the structure of the triangular pyramid, we trivially accept that continuous function $\alpha,\beta, t$ in respect belong to hom homomorphism f, F, T.
	\end{remark} 
	
	\begin{proof}
		According to the use of Baer criterion in the three dimensional pyramid we have 
		\[hom(R_G,E) \enskip o \enskip  hom(I_g,R_G)= hom(I_g,E)\]	
		if \[hom(R_G,E)=F\]
		\[hom(I_g,R_G)=T\]
		\[hom(I_g,E)=f\]
		According to this commutative  relation $FT=f$ of the triangular pyramid  we can see that $\beta\in F$ and $t\in T$ and  $\alpha\in f \Longrightarrow \beta t =\alpha$
	\end{proof}	
\end{theorem}

\begin{theorem}{2}
	In regards with topological space $T_b$ if   $I_g:C:\rightarrow \mathbb{R}$ is continuous function from a closed subset of C To $\mathbb{R}$  then there exists a function 
	
	\[E:T_b\rightarrow \mathbb{R} \quad  \text{such that}\quad E|_C= I_g \quad or \quad E\alpha =I_g\]
	
	\begin{center}

		\begin{tikzcd}
		C \arrow[r, "\alpha"] \arrow[rd, "I_g"'] & T_b \arrow[d, "E"] \\
		& \mathbb{R}        
		\end{tikzcd}
	\end{center}
	\begin{proof}
		According to the use of the Tietze extension theorem  in the three dimensional pyramid we have 
		
		\[(a )\qquad R_Gt=I_g  \]
		
		and from previous theorem we have 
		
		\[\ (b)\qquad \beta t=\alpha\]	
		now we proof this theorem
		
		\[\ (c) \qquad E\alpha =I_g\]
		
		We place $(a)$ and $(b)$ in $(c)$  and we obtain this relation
		\[(d) \qquad E\beta t=R_G t\]
		
		Based on the three dimensional triangular pyramid, The exiestence of (d) can be cleary proved.

	\end{proof}
\end{theorem}

\begin{theorem}
	In regards with topological space $T_b$ if $\mathbb{R}_G : N\rightarrow \mathbb{R}$ is a continuous function from normal topological space N to $\mathbb{R}$ then there exist a map
	\[E: T_b\rightarrow \mathbb{R}\quad \text{such that}\quad E|_N=R_G \quad or \quad  E\beta =R_G\]
	
	\begin{center}

		\begin{tikzcd}
		T_b \arrow[d, "E"'] & N \arrow[l, "\beta"'] \arrow[ld, "R_G"] \\
		\mathbb{R}          &                                        
		\end{tikzcd}
	\end{center}
	\begin{proof}	
		Since we have these results from preovious theormes 
		
		\[(a) \qquad \beta t=\alpha \Longrightarrow \beta=\alpha t^{-1}\]
		\[(b)\qquad R_Gt=I_g \longrightarrow R_G=I_g t^{-1}\]
		
		By placing these results in the relation $E\beta=R_G$ we have 
		$E\alpha t^{-1}=I_g t^{-1}$. With regrad to the previous theorem $E\alpha=I_g$ the existence of the last relation is proved.

	\end{proof}	
\end{theorem}

\section{Refrence}

\end{document}